\documentclass[a4paper,12pt]{amsart}
\usepackage[cp1251]{inputenc}
\usepackage[T2A]{fontenc}
\usepackage[russian,ukrainian,english]{babel}
\usepackage{amscd}
\usepackage{amssymb}
\usepackage{amsmath,amsthm,amsfonts}
\usepackage{hhline}

\usepackage{cite}

\newtheorem{theorem}{Theorem}[section]
\newtheorem{lem}{Lemma}[section]
\newtheorem{ozn}{Definition}[section]
\newtheorem{remk}{Remark}[section]
\newtheorem{nas}{Corollary}[section]
\newtheorem{prop}{Proposition}[section]

\newcommand{\rr}{\mathbb{R}}

\renewcommand{\emptyset}{\varnothing}
\newcommand{\Cl}[1]{\overline{#1}}
\newcommand{\Int}[1]{\mathrm{Int}\,{#1}}
\newcommand{\Fr}[1]{\mathrm{Fr}\,{#1}}

\newcommand{\restrict}[2]{#1\raisebox{-0.65ex}{$\left|\vphantom{{#1}_{#2}}\right.$}_{#2}}

\newcommand{\Rfr}{R_{\mathrm{Fr}}}

\usepackage{graphicx}

\begin{document}

\begin{center}
{\large\bf On the conditions of topological equivalence of pseudoharmonic functions defined on disk}
\end{center}
\vspace{3mm}
\begin{center}
{\small Polulyakh E.,\;Yurchuk I.}
\end{center}
\begin{center}
{\small Institute of Mathematics of Ukrainian National Academy of Sciences, Kyiv}
\end{center}
\vspace{5mm}

{\noindent \small\textbf{Abstract.} Let $D^{2}\subset \mathbb{C}$ be a closed two-dimensional disk and $f:D^2\rightarrow \mathbb{R}$ be a continuous function such that
 a restriction of $f$ to $\partial D^{2}$ is a continuous function with a finite number of local extrema and $f$ has a finite number of critical points in $\Int{D^{2}}$ such that each of them is saddle (i.e., in its neighborhood the local representation of $f$ is $f = Re z^n + const$, where $z=x+iy$, $n \geq 2$). This class of functions coincides with class of pseudoharmonic functions defined on $D^2$~\cite{Mr, st, Bo, MrJ, Kp, Mrt, BoM}.

First, we will construct an invariant of such functions which contains all information about them. Then, in terms of such invariant the necessary and sufficient conditions for pseudoharmonic functions to be topologically equivalent will be obtained.}
\medskip

{\noindent \small\textbf{Keywords}. a pseudoharmonic function, a combinatorial diagram, a topological conjugancy.}

\section*{Introduction}

In~\cite{arn, Ms, Os, Pr1, Sh1, Sh, ArD, SOB, BF, Pr2} the problems of topological classification of functions, vector fields and others structures on manifold were solved. In most cases, such solutions were received by the construction of combinatorial objects which contain all necessary facts about the structure being investigated. For example, in~\cite {SOB} authors constructed the spin graphs in order to classify Morse-Smale's fields on closed two-dimensional manifolds. So, an isomorphism of such graphs is the necessary and sufficient condition of topological classification of fields. In~\cite {arn, ArD} the classification of $M$--morsifications and bifurcations was obtained in terms of snakes (the special permutations).

Our goal is to construct an invariant of pseudoharmonic functions defined on disk and to receive the conditions for them to be topologically equivalent.

In Section 1 we introduce all necessary definitions which connect the nature of pseudoharmonic functions and theory of graphs.

In Section 2 the invariant of such functions is constructed, and also its main properties are studied. In particular, the invariant is a finite connected graph with a partial orientation and a partial order on its vertices which is generated by a function.

The main result of this paper is Theorem 3.1. in Section 3 which formulates the necessary and sufficient conditions for pseudoharmonic functions to be topologically equivalent in terms of their combinatorial diagrams.

\section{Preliminaries}

Let us remind some definitions and results which will be helpful for us.

Let $T$ be a tree with a set of vertices $V$ and a set of edges $E$. Suppose that $T$ is non degenerated ( has at least one edge). Denote by $V_{ter}$ a set of all terminal vertices of $T$ (i.e. vertices of degree one). Suppose that for some subset $V^{\ast} \subseteq V$ the following  condition holds true
\begin{equation}\label{eq_01}
V_{ter} \subseteq V^{\ast} \,.
\end{equation}
Let also $\varphi : T \rightarrow \rr^{2}$ is an embedding such that
\begin{equation}\label{eq_02}
\varphi(T) \subseteq D^{2} \,, \quad \varphi(T) \cap \partial D^{2} = \varphi(V^{\ast}) \,.
\end{equation}

\begin{lem}[see \cite{PY}]\label{lemma_01}
A set $\rr^{2} \setminus (\varphi(T) \cup \partial D^{2})$ has a finite number of connected components
\[
U_{0} = \rr^{2} \setminus D^{2}, U_{1}, \ldots, U_{m} \,,
\]
and for every $i \in \{1, \ldots, m\}$ a set $U_{i}$ is an open disk and is bounded by a simple closed curve
\[
\partial U_{i} = L_{i} \cup \varphi(P(v_{i}, v_{i}')) \,, \quad
L_{i} \cap \varphi(P(v_{i}, v_{i}')) = \{ \varphi(v_{i}), \varphi(v_{i}') \}
\]
 where $L_{i}$ is an arc of $\partial D^{2}$ such that the vertices $\varphi(v_{i})$ and $\varphi(v_{i}')$ are its endpoints, and  $\varphi(P(v_{i}, v_{i}'))$ is an image of the unique path $P(v_{i}, v_{i}')$ in $T$ which connects $v_{i}$ and $v_{i}'$.
\end{lem}

It is known~\cite{Z-F-C} that if $E_{1}$, $E_{2}$ are closed disks and $h:\partial E_{1}\rightarrow \partial E_{2} $ is a homeomorphism that there exists a homeomorphism $H:E_{1}\rightarrow E_{2}$ such that $H|_{\partial E_{1}}=h$.

\begin{ozn}
Two functions $f,g:D^{2}\rightarrow \mathbb{R}$ are called \emph{topologically equivalent} if there exist orientation preserving homeomorphisms $h_{1}:D^{2}\rightarrow  D^{2}$ and $h_{2}:\mathbb{R}\rightarrow \mathbb{R}$ such that $f=h^{-1}_{2}\circ g\circ h_{1}$.
\end{ozn}

We remind that function $f(x,y)$ is harmonic at a point $(x_0,y_0)$ if $\frac{\partial^2 f}{\partial x^2}(x_0,y_0) + \frac{\partial^2 f}{\partial y^2}(x_0,y_0)=0$.

\begin{ozn}
Function $f(z)$ is pseudoharmonic at a point $z_0=(x_0,y_0)$ if there exist a neighborhood $U(z_0)$ and a homeomorphism $\varphi$ of $U(z_0)$ onto itself such that $\varphi(z_0)=z_0$ and $f(\varphi(z))$, $z=(x,y)$, is harmonic.
\end{ozn}

Function $f$ is pseudoharmonic in a domain if it is pseudoharmonic at any point of it.

\begin{ozn}\label{ozn_f_1}
Point $z_{0} \in D^{2}$  is \emph{a regular} point of $f$ if there exist its open neighborhood $U \subseteq D^{2}$ and a homeomorphism $\varphi : U \rightarrow \Int{D^{2}}$ such that $\varphi(z_{0}) = 0$ and $f \circ \varphi^{-1}(z) = Re z + f(z_{0})$ for all $z \in U$.

The neighborhood $U$ will be called \emph{canonical}.
\end{ozn}

\begin{ozn}\label{ozn_f_2}
Point $z_{0} \in \partial D^2$ is \emph{a regular boundary point} of $f$ if there exist its neighborhood $U$ in $D^{2}$ and a homeomorphism $h : U \rightarrow D^{2}_{+}$ of such neighborhood into upper half-disk $D^{2}_{+}$ such that $h(z_{0}) = 0$, $h(U \cap f^{-1}(f(z_{0}))) = \{0\} \times [0, 1)$, $h(U \cap \partial D^2) = (-1, 1) \times \{0\}$ and a function $f \circ h^{-1}$ is strictly monotone on the interval $(-1, 1) \times \{0\}$.

The neighborhood $U$ will be named \emph{canonical}.
\end{ozn}

\begin{remk}\label{remk_kanon_neighb}
It is easy to see that the canonical neighborhoods from Definitions~\ref{ozn_f_1} and~\ref{ozn_f_2} can be chosen as small as need.
\end{remk}

If a point $z_{0} \in \Int D^{2}$ is not a regular point of $f \in F(D^{2})$ it will be called \emph{critical}.
By definition all critical points of $f$ are saddle.

Point of $\partial D^{2}$ that is neither a boundary regular point nor an isolated point of its level curve will be called \emph{a critical boundary point}.

\begin{ozn}
Number $c$ is \emph{a critical value of $f$} if level set $ f^{-1}(c)$ contain critical points.

Number $c$ is \emph{a regular value of $f$} if a level set $ f^{-1}(c)$ does not contain critical points and it is homeomorphic to a disjoint union of segments which intersect with a boundary $\partial D^{2}$ only in their endpoints.
\end{ozn}

It is known that any level curve of pseudoharmonic function is homeomorphic to a disjoint union of trees~\cite{Mr, Bo}.

\begin{ozn}
Number $c$ is \emph{a semiregular} value of $f$ if it is neither regular nor critical.
\end{ozn}

\begin{remk}
From Definitions it follows that level curves of semiregular value contain only boundary critical points and local extrema of $f$ (they belong to $\partial D^2$ and are isolated points of level curves of $f$). The level curves of the critical value contain the critical points and they also can contain boundary critical points and local extrema.
\end{remk}

From Theorem 4.1~\cite{Mrt}, see also~\cite{Mr}, it follows that for any critical boundary point there exists an homeomorphism of its canonical neighborhood onto half-disk which maps that point to origin and an image of its level set consists of finite number of rays outgoing from it.

It is easy to prove that the level curves of a semiregular value of pseudoharmonic function are isomorphic to trees, in general disconnected.

\begin{figure}[htbp]
\centerline{\includegraphics{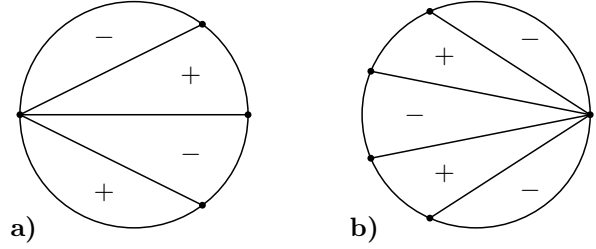}}
\label{fig2} \caption{ In case $a$) point is regular but in case $b$) it is local maximum of $f|_{\partial D^{2}}$.}
\end{figure}

For any boundary critical point $z_{0}$ the number of domains on which half-disk is divided by a set $f^{-1}(f(z_{0}))$ is greater than 2. Every such domain can be associated with a sign either ``$+$'' or ``$-$'' depending on a sign of the difference $f(z) - f(z_{0})$ in it.

We should remark that if a number of the domains is even then the domains adjoining to $\partial D^{2}$ have the same sign. Therefore a point $z_{0}$ is local extremum of $f|_{\partial D^{2}}$. It is obvious that in the case when this sign is minus it is a local maximum, otherwise it is minimum.

Let $n>1$ and $z_1,\ldots,z_{2n}$ is a sequence of points on $\partial D^2\cong S^1$ and they are passed in this order in the positive direction on $\partial D^2$. Let $\gamma_k$ is a positively oriented close arc of $\partial D^2 $ from $z_k$ to $z_{k+1}$ or $z_1$ when $k=2n$. Denote by
$\mathring{\gamma}_{k}$ an open arc $\gamma_k$ without its endpoints.

\begin{ozn}[see~\cite{Pol3}]\label{ozn_weak_regular}
Suppose that for a continuous function $f : D \rightarrow \rr$ there exists $n = \mathcal{N}(f) \geq 2$ and a sequence of points $z_1, \ldots, z_{2n-1}, z_{2n} \in \Fr{D}$ (in this order these points are passed in the positive direction on $\Fr{D}$) such that the following conditions are satisfied:
\begin{itemize}
	\item [1)] every point of the domain $\Int{D} = D \setminus \Fr{D}$ is a regular point of $f$;
	\item [2)] $\mathring{\gamma}_{2k-1} \neq \emptyset$ for $k \in \{1, \ldots, n\}$, and any point of the arc $\mathring{\gamma}_{2k-1}$ is a regular boundary point of $f$ (in particular, the restriction of $f$ to $\gamma_{2k-1}$ is strictly monotone);
	\item [3)] the arcs $\gamma_{2k}$, $k \in \{1, \ldots, n\}$, are connected components of the level sets of $f$.
\end{itemize}
Such function will be called \emph{a weakly regular on $D$}.
\end{ozn}

\begin{lem}[see~\cite{Pol3}]\label{lem_number_of_arks}
Let $f$ be a weakly regular function on $D$. Then $\mathcal{N}(f) = 2$.
\end{lem}

\begin{ozn}[see~\cite{Pol3}]\label{ozn_almost_weak_regular}
Suppose that for some $n \geq 2$ and a sequence of points $z_1, \ldots, z_{2n} \in \Fr{D}$ a function $f$ satisfies all conditions of Definition~\ref{ozn_weak_regular} except 3, but instead, the following condition holds true
\begin{itemize}
	\item[$3')$] for $j = 2k$, $k \in \{1, \ldots, n\}$, an arc $\gamma_j$ belongs to some level set of $f$.
\end{itemize}
Such function will be named \emph{ an almost weakly regular on $D$}.
\end{ozn}

Let $f$ be an almost weakly regular function on $D$. Denote by $2\cdot\mathcal{N}(f)$ a minimal number of points and arcs satisfying Definition~\ref{ozn_almost_weak_regular}. Obviously, $2\cdot\mathcal{N}(f)$ only depends on $f$.

\begin{prop}[see~\cite{Pol3}]\label{prop_min_num_of_arcs}
If for some $n \geq 2$ and a sequence of points $z_1, \ldots, z_{2n} \in \Fr{D}$ a function $f$ satisfies all conditions of Definition~\ref{ozn_almost_weak_regular} and $n = \mathcal{N}(f)$, then a collection of sets $\{\mathring{\gamma}_{2k-1}\}_{k=1}^n$ coincides with a collection of connected components of the set of regular boundary points of $f$.
\end{prop}

\begin{ozn}[see~\cite{Pol3}]
A simple continuous curve $\gamma:[0,1]\rightarrow \overline{W}$ is called an \emph{$U$-trajectory} if $f\circ \gamma$ is strongly monotone on the segment $[0,1]$.
\end{ozn}

\begin{ozn}[see~\cite{Pol3}]
Let $f$ be a weakly regular function on the disk $D$, let $\gamma_1,\ldots,\gamma_4$ be arcs from Definition~\ref{ozn_weak_regular}. If through every point of a set $\Gamma$ which is dense in $\mathring{\gamma}_{2}\cup \mathring{\gamma}_{4}$ passes a $U$-trajectory, then the function $f$ is called \emph{regular on} $D$.
\end{ozn}

\begin{theorem}[see~\cite{Pol3}]\label{theorem_1.3}
Let $f$ be a regular function on $D$, $\gamma_1, \ldots, \gamma_4$ be the arcs from Definition~\ref{ozn_weak_regular}.

Let $D' = I^2$, if $\mathring{\gamma}_2 \neq \emptyset$ and $\mathring{\gamma}_4 \neq \emptyset$; $D' = D^2$, if $\mathring{\gamma}_2 \cup \mathring{\gamma}_4 = \emptyset$; $D' = \Cl{D}^{2}_{+}$, if exactly one of sets either $\mathring{\gamma}_2$ or $\mathring{\gamma}_4$ is empty.

Suppose that $\phi : \Fr{D} \rightarrow \Fr{D'}$ is a homeomorphism such that $\phi(K) = K'$, where
\begin{gather*}
K = f^{-1}\bigl(\min_{z \in D}(f(z)) \cup \max_{z \in D}(f(z))\bigr) \,,\\
K' = \Bigl\{ (x, y) \in D' \,\bigl|\, y \in  \bigl\{\min_{(x, y) \in D'}(y), \max_{(x, y) \in D'}(y) \bigr\}\bigr.\Bigr\} \,.
\end{gather*}

There exists a homeomorphism $H_f$ of disk $D$ onto $D'$ such that $H_f|_K = \phi$ and $f \circ H_f^{-1}(x, y) = ay + b$, $(x, y) \in D'$, for some $a, b \in \rr$, $a \neq 0$.
\end{theorem}


\section{Combinatorial invariant of pseudoharmonic functions}

At first we should remind the term of Reeb's graph.
Let $M$ be a smooth compact manifold. Suppose that $f:M\rightarrow\mathbb{R}$ is a smooth function with a finite number of critical points.
Let us define connected component of level curves of $f^{-1}(a)$, where $a\in\mathbb{R}$, as layer. Then $ M$ is the union of all layers of $f$.
Also we can define the relation of equivalence as the property of points to belong to a same layer and consider the quotient space by this relation. It is homeomorphic to a finite graph named Reeb's graph and let us denote it by $\Gamma_{K-R}(f)$. Its vertices are components of level curves such that
they contain the critical points.

Let $D^2$ be a closed oriented disk and $f:D^2\rightarrow \mathbb{R}$ be a pseudoharmonic function. We should remark that for a manifold with boundary the construction of Reeb's graph is an open problem therefore there is a reason to obtain another
invariant for such functions.

\textbf{\textit{Construction of invariant}} for pseudoharmonic function named as combinatorial diagram:
\begin{enumerate}
	\item[1)] We construct Reeb's graph $\Gamma_{K-R}(f|_{\partial D^{2}})$ of the restriction of $f$ to $\partial D^{2}$. It is isomorphic to
 circle with even number of vertices of degree 2 (vertices are local extrema of the restriction of $f$ to $\partial D^{2}$) and fix an orientation on $\Gamma_{K-R}(f|_{\partial D^{2}})$ which is generated by the orientation of $D^2$.
	\item[2)] Let $a_{i}$ be the critical values of $f$ and $c_j$ be the semiregular values. We add to $\Gamma_{K-R}(f|_{\partial D^{2}})$ those
connected components of sets
	\[
	f^{-1}(a_{1})\bigcup \ldots \bigcup f^{-1}(a_{k})\bigcup f^{-1}(c_1)\bigcup f^{-1}(c_2)\bigcup\ldots\bigcup f^{-1}(c_l) \,,
	\]
	of level curves that contain critical and boundary critical points. It is obvious that new vertices appear on $\Gamma_{K-R}(f|_{\partial D^{2}})$.
We set
	\[
	P(f)=\Gamma_{K-R}(f|_{\partial D^{2}}) \cup \bigcup\limits_{i}\widehat{f}^{-1}(a_{i}) \cup \bigcup\limits_{j}\widehat{f}^{-1}(c_j) \,,
	\]
	where $\widehat{f}^{-1}(a_{i})\subset f^{-1}(a_{i})$, $\widehat{f}^{-1}(c_j)\subset f^{-1}(c_j)$ are those connected components of level sets
that contain critical and boundary critical points.
	\item[3)] We put a partial order on vertices of $P(f)$ by using the values of $f$: $v_{1} < v_{2}$ $\Longleftrightarrow$ $f(x_{1}) < f(x_{2})$,
where $v_{1},v_{2}\in P(f)$, $x_{1},x_{2}$ are points corresponding to vertices $v_{1},v_{2}$, respectively. In case of the same values of function
on vertices they will be non comparable.
\end{enumerate}
This partial order is strict~\cite{Mel} since the relation is antireflexive, antisymmetric and transitive.
$P(f)$ will be called \emph{combinatorial diagram} of pseudoharmonic function $f$.

By the construction  $P(f)$ is a finite partially oriented graph with a strict partial order on vertices.
\begin{figure}[htbp]
\begin{center}
\includegraphics{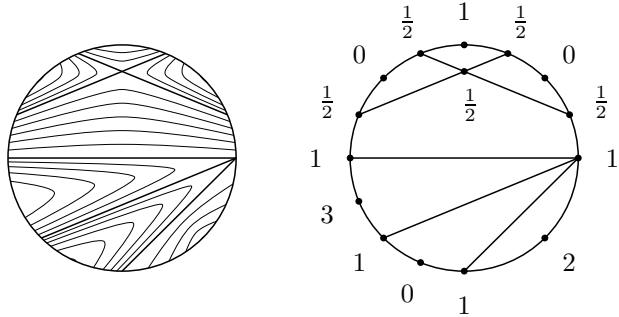}
\caption{Example of a diagram of some pseudoharmonic function.}
\end{center}
\end{figure}

We constructed the combinatorial invariant of $f$ as subset of $D^2$. We will consider it as the abstract partially oriented graph with fixed relation of partial
order on the set of vertices $V(P(f))$.
\begin{ozn}
Two combinatorial diagrams $P(f)$ and $P(g)$ are \emph{isomorphic} if there exists an isomorphism $\phi : P(f) \rightarrow P(g)$ between
 them which preserves a strict partial order given on their vertices (maps $\restrict{\phi}{V(P(f))}$ and
  $\phi^{-1}\restrict{\vphantom{\phi}}{V(P(g))}$ are monotone) and the orientation.
\end{ozn}
We put the natural topology on the diagram $P(f)$. For example, it can be introduced by structure of one-dimensional CW-complex on $P(f)$.
All vertices of the graph $P(f)$ can be considered as 0-dimensional cells, similarly, all edges can be considered as 1-dimensional cells.
$P(f)$ also can be regarded as a subset of $\rr^3$  and all edges are straight segments.
\begin{ozn}

Homeomorphism $\varphi: P(f) \rightarrow P(g)$ is said \emph{to realize an isomorphism $\phi: P(f) \rightarrow P(g)$ of combinatorial diagrams}
 if $\restrict{\varphi}{V(P(f))} = \restrict{\phi}{V(P(f))}$ and from $\phi(e) = e'$ it follows that $\varphi(e) = e'$ for any edge $e \in E(P(f))$.
\end{ozn}
\begin{remk}\label{remk_non_unique_realiz}

It is clear that every isomorphism $\phi$ of combinatorial diagrams is realized by some homeomorphism but it is not uniquely defined:
for every edge $e \in E(P(f))$ we can arbitrarily choose a homeomorphism $\varphi_e : e \rightarrow \phi(e)$ such that maps $\varphi_e$
and $\phi$ are the same on $e \cap V(P(f))$.
\end{remk}

We constructed the combinatorial diagram $P(f)$ as a subset of $D^2$ therefore ``support'' of diagram in $D^2$ is correctly defined since it is the set
\begin{equation}\label{eq_bearer}
P_f = \Gamma_{K-R}(f|_{\partial D^{2}}) \cup \bigcup\limits_{i}\widehat{f}^{-1}(a_{i}) \cup \bigcup\limits_{j}\widehat{f}^{-1}(c_j) \,,
\end{equation}
 where $\widehat{f}^{-1}(a_{i})$ and $\widehat{f}^{-1}(c_j)$ are connected components of level curves of $f$ which contain the critical and
 boundary critical points.

Similarly, to the vertices of $P(f)$ corresponds the set $V_f$ which is the ``support'' of the set of its vertices in $D^2$. Function $f$ induces
a strict partial order on it using the following correlations  $x_1 < x_2 \Leftrightarrow f(x_1) < f(x_2)$.
Denote by $M(f) \subset \partial D^2$ \emph{the set of local extrema of $f$ on $D^2$}. By the construction every point of this set corresponds
to some vertex of $P(f)$, thus $M(f) \subset V_f$. Other vertices of $P(f)$ are characterized by  the property that each of them is a common endpoint of at least three edges, therefore it has no neighborhood that is homeomorphic to segment in the space $P(f)$.

\begin{ozn}\emph{$\mathcal{C}r$-subgraph} of $P(f)$ is a subgraph $q(f)$ such that:
\begin{itemize}
\item $q(f)$ is a simple oriented cycle;
\item arbitrary pair of adjacent vertices $v_{i},v_{i+1}\in q(f)$ is comparable.
\end{itemize}
\end{ozn}

Let $\varphi : P(f) \rightarrow D^2$ be an arbitrary embedding of topological space $P(f)$ into $D^2$ such that $\varphi(P(f)) = P_f$. Granting
what we said above it is obvious that an inclusion $M(f) \subseteq \varphi(V(P(f)))$ is equivalent to $\varphi(V(P(f))) = V_f$.

\textit{\textbf{In what follows unless otherwise stipulated we assume that for any embedding of $P(f)$ into $D^2$ the orientation of $\mathcal{C}r$-subgraph coincides with the orientation of $\partial D^2$.}}

\begin{ozn}
Let $\varphi : P(f) \rightarrow D^2$ be an embedding of topological space $P(f)$ into $D^2$. It is called \emph{to be consistent with $f$} if
 the following correlations hold true:
\begin{itemize}
	\item $\varphi(P(f)) = P_f$;
	\item $M(f) \subseteq \varphi(V(P(f)))$;
	\item a partial order on $(V(P(f))) = V_f$ induced  by a partial order on $V(P(f))$ with help of $\varphi$ coincides with a partial
order induced on this set from $\rr$ by $f$.
\end{itemize}
\end{ozn}

It is clear that there exist at least one embedding $\varphi : P(f) \rightarrow D^2$ which is consistent with $f$. If $\psi : P(f) \rightarrow P(f)$
is an isomorphism of $P(f)$ onto itself (for example, identical map) which can be realized by homeomorphism $\hat{\psi} : P(f) \rightarrow P(f)$,
then an embedding $\varphi \circ \hat{\psi}$ is also consistent with $f$.

We should remind that vertices $v_{1}$ and $v_{2}$  of some graph $G$ are \emph{adjacent} if they are endpoints of the same edge.

Let $v$ be some vertex of the diagram $P(f)$ and $\{v_{i}\}$,  $i=\overline{1,k}$, be a set of all adjacent vertices to it. Then there exist
points $x$ and $x_{i}$ of $D^{2}$ that correspond to vertices $v$ and $v_{i}$. Denote by $X_{i}\subseteq D^2$  the set of points which corresponds
to edge $e(v,v_{i})$ (it is clear that every $X_{i}$ is homeomorphic to segment). Let us consider the following cases:

\emph{Case 1: $x\in \Int D^{2}$}.
 Then $f(x)=f(x_{i})=a$, where $i=\overline{1,k}$ and $a$ is a critical value. Therefore vertices $v$,$v_{1},v_{2},\ldots,v_{k}$ are pairwise non
 comparable. Since level set of the critical value $a$ is a finite tree then all vertices of it are non comparable.

\emph{Case 2: $x\in \partial D^{2}$}.
In this case the point $x$ is either regular or local extremum of $f|_{\partial D^{2}}$ which is continuous and monotonically increase
(decrease) between adjacent local extrema. Therefore, among sets $X_{i}$ there exist such that function monotonically increases (decreases)
on them. Circle is closed Jordan curve then there are exactly two such sets $X_{j}$ and $X_{k}$ whose endpoints are points $x_{j}$ and $x_{k}$.
So, it follows that among all vertices $\{v_{i}\}$ adjacent to $v$ there exist exactly two vertices $v_{j}$ and $v_{k}$ which are comparable with a vertex $v$.
For both $v_{j}$ and $v_{k}$ there exist exactly two vertices which are comparable to it thus these vertices generate a cycle (the case of two  or more non intersecting cycles is impossible since a disk has one boundary circle).

It is obvious that $v$ together with both vertices $v_{j}$ and $v_{k}$ belong to $q(f)$--cycle.

The fact that the diagram $P(f)$ is constructed by pseudoharmonic function implies several characteristics of it.

\textbf{\textit{Main properties of $P(f)$}}:
\begin{enumerate}
	\item[C1)]\label{c1} there exists the unique $\mathcal{C}r$-subgraph $q(f)\in P(f) $;
	\item[C2)] $\overline{P(f)\setminus q(f)}=\bigcup\limits_{i} \Psi_{i}$, $\Psi_{j}\bigcap\Psi_{i}=\varnothing$, where $i\neq j$,
and every $\Psi_{i}$ is a tree such that for any index $i$ arbitrary two vertices $v',v''\in \Psi_{i}$ are non comparable;
	
	\item[C3)] there exists an embedding $\psi:P(f)\rightarrow D^{2}$ such that $\psi(q(f))=\partial D^{2}$ and
$ \psi(P(f)\setminus q(f))\subset \Int D^{2}$;
	\item[C4)] for every connected component $\Theta$ of $D^2 \setminus P_f$ the function $f$ is regular (see~\cite{Pol3}) on the set $\Cl{\Theta}$.
\end{enumerate}
From what was said above the existence of $\mathcal{C}r$-subgraph and the fairness of $C2$ follow. From the existence of
$\mathcal{C}r$-subgraph and $C2$ it follows that $q(f)$ is unique.
Condition $C3$ follows from fact that $P(f)$ is a diagram of a function $f$, defined on $D^{2}$. $\mathcal{C}r$-subgraph $q(f)\in P(f)$ is unique thus from the
definitions it is easy to see that for every embedding $\psi : P(f) \rightarrow D^2$ which is consistent with $f$ the
equality $\psi(q(f)) = \partial D^2$ should hold true.

By the definition of the diagram $P(f)$ any tree $\Psi_i$ corresponds to a connected component of some critical or semiregular
level set of $f$. A number of trees is the same as a number of such components which contain critical or boundary critical points.
Denote by $P_f^c = \Cl{P_f \setminus \partial D^2}$ the union of such components.

Let $\psi : P(f) \rightarrow D^2$ be an embedding which is consistent with $f$.
If the endpoints $v'$ and $v''$ of some edge $e = e(v', v'')$ of $P(f)$ are non comparable, then
$\mathring{e} = e \setminus \{v', v''\} \in P(f) \setminus q(f) \subseteq \bigcup_i \Psi_i$. Thus
$\psi(\mathring{e}) \subseteq P_f \cap \Int{D^2} \subseteq P_f^c$. Then there exists
$c = c(e) \in \rr$ such that $\varphi(e) \subset f^{-1}(c)$. Any connected set $\psi(\Psi_i)$ belongs to some connected component of $P_f^c$.
From the facts that a map $\psi$ is an embedding and $\psi(q(f)) = \partial D^2$ follow the equalities
\[
\psi\Bigl(\bigcup_i \Psi_i\Bigr) = \psi\bigl(\,\Cl{P(f) \setminus q(f)}\,\bigr) = \Cl{\psi(P(f)) \setminus \psi(q(f))} = \Cl{P_f \setminus \partial D^2} = P_f^c \,.
\]
By the definition the number of connected components of sets $\bigcup_i \Psi_i$ and $P_f^c$ coincides  thus any set $\psi(\Psi_i)$ is a connected component of $P_f^c$.

Let us combine together corollaries of Conditions $C1$--$C3$ which we obtained above.
\begin{prop}\label{prop_2.1}
Let $P(f)$ be a combinatorial diagram of pseudoharmonic function $f$ and $\psi : P(f) \rightarrow D^2$ be an embedding which is consistent with $f$.
Then the following conditions hold true:
\begin{itemize}
	\item $\psi(q(f)) = \partial D^2$;
	\item for any tree $\Psi_i$ the set $\psi(\Psi_i)$ is a component of critical or semiregular level set of $f$.
\end{itemize}
\end{prop}

Let us prove Condition $C4$.
\begin{prop}\label{prop_2.2}
Let $P(f)$ be a combinatorial diagram of pseudoharmonic function $f$ and $\psi : P(f) \rightarrow D^2$ be an embedding such that $\psi(q(f)) = \partial D^2$.

The set $\Fr{\Sigma} = \Fr{\Cl{\Sigma}}$ is an image of a simple cycle $Q$ of $P(f)$ for any connected component $\Sigma$ of $D^2 \setminus \psi(P(f))$.
\end{prop}
\begin{proof}
All vertices of $\Psi_j \subset P(f)$ which do not belong to $\mathcal{C}r$-cycle $q(f)$ correspond to critical points of $f$ for any $j$, thus they have even degree no smaller than 2. Therefore the set $V_{ter}^j$ of all vertices of $\Psi_j$ of degree 1 is contained in $q(f)$ and we can apply Lemma~\ref{lemma_01} to a map $\restrict{\psi}{\Psi_j}$.

By induction on the number of trees $\Psi_i$ embedded into disk from Lemma~\ref{lemma_01} it follows that a boundary $\Fr{\Sigma}$ of $\Sigma$ is simple Jordan curve. Let us prove that its preimage $Q = \psi^{-1}(\Fr{\Sigma})$ is a subgraph of $P(f)$.
It suffices to verify the following assertion. Let $e = e(v_1, v_2)$ be some edge of $P(f)$ and $x \in \mathring{e} = e \setminus \{v_1, v_2\}$ be an inner point of $e$. If $x \in Q$, then $e \subset Q$.

It is obvious that the set $Q$ is a simple closed curve. Therefore $Q \setminus \{x'\}$ is connected for any $x' \in Q$. Thus $Q \setminus e \neq \emptyset$ (any point of segment $e$ except its endpoints splits it, see~\cite{Newman}).
Suppose that an edge $e$ is support of simple continuous curve $\alpha : I \rightarrow P(f)$, $\alpha(0) = v_1$, $\alpha(1) = v_2$. Therefore $x = \alpha(\tau)$ for some $\tau \in (0, 1)$. We should remark that $e$ is one-dimensional cell of CW-complex $P(f)$, thus $\mathring{e} = \alpha(\mathring{I})$ is an open subset of $P(f)$ (we denoted $\mathring{I} = (0, 1)$). For any interval $\mathring{I}(t_1, t_2) = (t_1, t_2)$, $t_1, t_2 \in I$, $t_1 < t_2$, the set $\alpha(\mathring{I}(t_1, t_2))$ is an open subset of $P(f)$. It is also obvious that $P(f) \setminus \alpha([t_1, t_2])$, where $t_1, t_2 \in I$, $t_1 < t_2$ is open in $P(f)$.

Let us show that at least one of sets $\alpha([0, \tau])$, $\alpha([\tau, 1])$ belong to $Q$. Suppose that it does not hold true. So, there exist $t_1 \in [0, \tau]$ and $t_2 \in [\tau, 1]$ such that $\alpha(t_1), \alpha(t_2) \notin Q$. Then the nonempty  sets $Q \cap \alpha(\mathring{I}(t_1, t_2)) \ni x$ and $Q \setminus \alpha([t_1, t_2]) \supseteq Q \setminus e$ open in subspace $Q$ of  $P(f)$ generate a partition of $Q$, but it is impossible since $Q$ is connected.

Suppose that $\alpha(t) \notin Q$ for some $t \in I$. Without loss of generality we can assume that $t < \tau$. Then $\alpha([\tau, 1]) \subset Q$. Let us fix $t' \in (\tau, 1)$ and set $x' = \alpha(t')$. The nonempty open in $Q$ sets $Q \cap \alpha(\mathring{I}(t, t')) \ni x$ and $Q \setminus \alpha([t, t']) \supset Q \setminus e$ generate in $Q$ the partition of subset $Q \setminus \{x'\}$, but it is impossible since $Q \setminus \{x'\}$ is connected.

Thus $e \in Q$ and $Q$ is a subgraph of $P(f)$. The set $Q$ is homeomorphic to circle thus it is a simple cycle.
\end{proof}

\begin{lem}\label{lemma_2_C4}
Let $P(f)$ be a diagram constructed at pseudoharmonic function $f$ and $\psi : P(f) \rightarrow D^{2}$ be an embedding that is consistent with $f$.

Then for any component $\Theta$ of the complement $D^2 \setminus P(f) = D^2 \setminus P_f$ its closure $\Cl{\Theta}$ is homeomorphic to disk and $f$ is regular on $\Cl{\Theta}$.
\end{lem}
\begin{proof}

Let $\Theta$ be a connected component of $D^2 \setminus P(f)$. Let us prove, at first, that $f$ is weekly regular in $\Cl{\Theta}$.
From Propositions~\ref{prop_2.1} and~\ref{prop_2.2} it follows that a boundary of $\Theta$ is a simple closed curve. Thus $\Cl{\Theta}$ is a closed disk and
\[
\Fr{\Cl{\Theta}} = \Cl{\Theta} \cap P_f = (\Cl{\Theta} \cap \psi(q(f))) \cup \Bigl(\Cl{\Theta} \cap \psi\Bigl(\bigcup_i \Psi_i\Bigr)\Bigr) = \Gamma_V \cup \Gamma_E \cup \Gamma_T \,,
\]
where $\Gamma_T = \Cl{\Theta} \cap \psi(\bigcup_i \Psi_i)$; $\Gamma_V = \Cl{\Theta} \cap M(f)$ is a set of points of $\partial{D^2} \cap \Fr{\Cl{\Theta}}$ which correspond to vertices of $P(f)$ of $q(f) \setminus \bigcup_i \Psi_i$; $\Gamma_E$ are the open arcs of $\partial{D^2} \cap \Fr{\Cl{\Theta}}$ which correspond to the edges of the cycle $q(f)$ without endpoints. It is obvious that the sets $\Gamma_V$, $\Gamma_E$ and $\Gamma_T$ are pairwise disjoint.

The set $\Gamma_V$ consists of the isolated points of level sets of $f$. Each of them is a local extremum of $f$ in $D^2$.
The function $f$ is locally constant on $\Gamma_T$ therefore any connected component $K$ of such set belongs to $\psi\bigl(\bigcup_i \Psi_i\bigr)$ and there exists $c_K \in \rr$ such that $K \in f^{-1}(c_K)$. Let $\Gamma_K$ be a connected component of $f^{-1}(c_K) \cap \Cl{\Theta}$ containing $K$. Then $\Gamma_K \subseteq \Cl{\Theta} \cap \psi\bigl(\bigcup_i \Psi_i\bigr) = \Gamma_T$. Consequently $\Gamma_K = K$.

From the definition it follows that all points of $\Gamma_E$ are regular boundary points of $f$ in $D^2$. It is easy to see that sufficiently small canonical neighborhood of any point of $\Gamma_E$ belongs to $\Cl{\Theta}$ therefore all points of $\Gamma_E$ are regular boundary points of $f$ in $\Cl{\Theta}$.

By the definition the set $\Gamma_E$ has a finite number of connected components (their number is no more than a number of the edges of the cycle $q(f)$) therefore there exists a finite collection of points $z_1, \ldots, z_{2n} \in \Fr{\Cl{\Theta}}$ which divide the circle $\Fr{\Cl{\Theta}}$ into arcs $\gamma_1, \ldots, \gamma_{2n}$ such that $\Gamma_E = \bigcup_{k=1}^n \mathring{\gamma}_{2k-1}$ (some arcs with even indices can degenerate into points).

It is clear that $\Fr{\Cl{\Theta}} \setminus \bigcup_{k=1}^n \mathring{\gamma}_{2k-1} = \bigsqcup_{k=1}^n \gamma_{2k} = \Gamma_V \cup \Gamma_T$. The sets $\Gamma_V$ and $\Gamma_T$ are closed and disjoint thus any arc $\gamma_{2k}$, $k \in \{1, \ldots, n\}$, belongs to either $\Gamma_V$ or $\Gamma_T$.

From the preceding it follows that any set $\gamma_{2k}$, $k \in \{1, \ldots, n\}$, is a connected component of some level set of $f$ on $\Cl{\Theta}$. Therefore the collection of points $z_1, \ldots, z_{2n}$ satisfies to Definition~\ref{ozn_weak_regular} and $f$ is weakly regular on $\Cl{\Theta}$.

From Lemma~\ref{lem_number_of_arks} it follows that $n = \mathcal{N}(\restrict{f}{\Cl{\Theta}}) = 2$.
If $\mathring{\gamma}_{2k} \neq \emptyset$, $k \in \{1, 2\}$, then $\gamma_{2k} \in \Gamma_T$ (the set $\Gamma_V$ is discrete therefore $\gamma_{2k} \cap \Gamma_V = \emptyset$, see above) and any point $z \in \mathring{\gamma}_{2k}$ either belongs to $\Int{D^2}$ or is boundary critical point of $f$.

If $z \in \mathring{\gamma}_{2k} \cap \Int{D^2}$, then there exist an open neighborhood $W_z$ of $z$ in $D^2$ and a homeomorphism $\Phi_z : W_z \rightarrow \Int{D^2}$ such that $\Phi_z(z) = 0$ and $f \circ \Phi^{-1}_z(w) = Re w^m + f(z)$ for some $m \geq 2$.
The set $\Phi_z(f^{-1}(f(z)))$ divides $\Int{D^2}$ onto $2m$ open sectors such that each of them (for sufficiently small neighborhood $W_z$) belongs to $D^2 \setminus P(f)$. Thus for at least one of them its image under the action of $\Phi_z^{-1}$ belongs to $\Theta$. It is obvious that for every such sector there exists $U$-trajectory of $f$ which passes through the point $z$ and is contained in the closure of the image of sector under the action of $\Phi_z^{-1}$. Taking that into account some $U$-trajectory in $\Cl{\Theta}$ passes through $z$.

The number of the boundary critical points of $f$ on $D^2$ is finite therefore $\Gamma = \mathring{\gamma}_{2k} \cap \Int{D^2}$ is a dense subset of an arc $\mathring{\gamma}_{2k}$, $k \in \{1, 2\}$, and function $f$ is regular on $\Cl{\Theta}$.
\end{proof}

\begin{lem}\label{lemma_2_2}
Let $P(f)$ be a combinatorial diagram of pseudoharmonic function; $\psi_{1}$, $\psi_{2} : P(f) \rightarrow D^{2}$ be embeddings such that $\psi_{i}(q(f)) = \partial D^{2}$, $i = 1, 2$.

If an image $\psi_{1}(Q)$ of a simple cycle $Q \subset P(f)$ is a boundary of some component of the complement $D^{2} \setminus \psi_{1}(P(f))$, then an image $\psi_{2}(Q)$ is a boundary of some component of the complement $D^{2} \setminus \psi_{2}(P(f))$.
\end{lem}
\begin{proof}
Let us fix an embedding $\varphi : P(f) \rightarrow D^2$ consistent with $f$. Let $\psi : P(f) \rightarrow D^2$ be an embedding such that $\psi(q(f)) = \partial D^2$.
It is obvious that lemma follows from the following statement: an image $\psi(Q)$ of a simple cycle $Q \subseteq P(f)$ is a boundary of some component of the complement $D^2 \setminus \psi(P(f))$ iff a curve $\varphi(Q)$ is a boundary of one of components of $D^2 \setminus P_f = D^2 \setminus \varphi(P(f))$. Let us prove this statement.

Suppose that an image $\varphi(Q)$ of the cycle $Q$ bounds one of the components $\Theta$ of the set $D^2 \setminus P_f$. From Lemma~\ref{lemma_2_C4} it follows that $f$ is regular on disk $\Cl{\Theta}$, therefore there exist points $z_1, \ldots, z_4 \in \Fr{\Cl{\Theta}} = \varphi(Q)$ which divide a curve $\varphi(Q)$ into arcs $\gamma_1, \ldots, \gamma_4$ satisfying the following conditions:
\begin{itemize}
	\item $\mathring{\gamma}_1 \neq \emptyset$, $\mathring{\gamma}_3 \neq \emptyset$, and the set $\mathring{\gamma}_1 \cup \mathring{\gamma}_3$ is the set of boundary regular points of $f$ on $\Cl{\Theta}$;
	\item $\gamma_2$ and $\gamma_4$ are the components of level sets of $f$ on $\Cl{\Theta}$.
\end{itemize}
From these conditions it follows that (see a proof of Lemma~\ref{lemma_2_C4})
\begin{eqnarray}\label{eq_13}
\gamma_2 \cup \gamma_4 & = & \Fr{\Theta} \cap \Bigl( \varphi(V(P(f))) \cup \varphi \Bigl( \bigcup_i \Psi_i \Bigr) \Bigr) \supseteq \nonumber \\
	& \supseteq & \Fr{\Theta} \cap \varphi \Bigl( \bigcup_i \Psi_i \Bigr) = \Fr{\Theta} \cap \Bigl(\bigcup_i \varphi ( \Psi_i) \Bigr) \,.
\end{eqnarray}

Suppose that $\gamma_2 \subseteq f^{-1}(c')$, $\gamma_4 \subseteq f^{-1}(c'')$ for $c'$, $c'' \in \rr$. As $\mathring{\gamma}_1 \neq \emptyset$ and all points of this set are regular boundary points of $f$ on $\Cl{\Theta}$ the following statement holds true: $z_1 \neq z_2$ (since $\mathring{\gamma}_1 = \gamma_1 \setminus \{z_1, z_2\}$) and $f$ is strictly monotone on $\gamma_1$. Thus $c'' = f(z_1) \neq f(z_2) = c'$ and the sets $\gamma_2$ and $\gamma_4$ belong to different level sets of $f$.

Let us set $w_i = \psi \circ \varphi^{-1}(z_i)$, $\nu_i = \psi \circ \varphi^{-1}(\gamma_i)$, $i \in \{1, \ldots, 4\}$. The curve $\psi(Q)$ bounds an open domain $\Sigma$. From~\eqref{eq_13} it follows that $\nu_2 \cup \nu_4 \supseteq \Fr{\Sigma} \cap \psi (\bigcup_i \Psi_i)$.

Let us assume that a curve $\psi(Q) \subseteq \psi(P(f))$ is not a boundary of connected component of $D^2 \setminus \psi(P(f))$. Therefore $\Sigma \cap \psi(P(f)) \neq \emptyset$. We fix $z \in \Sigma \cap \psi(P(f))$. It is obvious that $\Sigma \subseteq \Int{D^2}$, therefore $x = \psi^{-1}(z) \in P(f) \setminus q(f) \subseteq \bigcup_i \Psi_i$ and $z \in \psi(\Psi_j)$ for some $j$.
All vertices of the tree $\Psi_j$ which do not belong to $\mathcal{C}r$-cycle $q(f)$ correspond to the critical points of $f$ thus they have even degree no less than 2. So, the set $V_{ter}^j$ of all vertices of degree one of tree $\Psi_j$ is contained in $q(f)$. It is easy to see that this set has at least two elements.

By easy check we can see that for any point $u$ of a subspace $\Psi_j$ of the space $P(f)$ there exist $v_u'$, $v_u'' \in V_{ter}^j$ and path $P(v_u', v_u'') \subset \Psi_j$ which connects $v_u'$ with $v_u''$ and passes through $u$.
Let us fix for a point $x = \psi^{-1}(z)$ vertices $v_x'$, $v_x'' \in V_{ter}^j$ and a path $P(v_x', v_x'') \subset \Psi_j$ which connects them and passes through a point $x$. We also fix a simple continuous curve $\alpha : I \rightarrow P(f)$ whose support is a path $P(v_x', v_x'')$. Suppose that $\alpha(0) = v_x'$, $\alpha(1) = v_x''$, $\alpha(\tau) = x$.

It is known that $\psi(V_{ter}^j) \subset \psi(q(f)) = \partial D^2$, but $z = \psi(x) \in \Sigma \subseteq \Int{D^2}$. Therefore $\tau \in (0, 1)$. Furthermore $\psi(v_x')$, $\psi(v_x'') \notin \Sigma$, so that each of the sets $\psi \circ \alpha([0, \tau])$ and $\psi \circ \alpha([\tau, 1])$ should intersect $\psi(Q) = \Fr{\Sigma}$. Suppose that
\begin{eqnarray*}
t' & = & \inf \{ t \in [0, \tau] \,|\, \psi \circ \alpha([t, \tau]) \in \Sigma \} \,,\\
t'' & = & \sup \{ t \in [\tau, 1] \,|\, \psi \circ \alpha([\tau, t]) \in \Sigma \} \,.
\end{eqnarray*}
Then $\tau \in (t', t'') \subseteq (\psi \circ \alpha)^{-1}(\Sigma)$ but $\psi \circ \alpha(t')$, $\psi \circ \alpha(t'') \in \psi(Q) = \Fr{\Sigma}$.
It is clear that for every $e = e(w', w'') \in P(v_x', v_x'')$ there exist $t', t'' \in I$, $t' < t''$ such that $w' = \alpha(t')$, $w'' = \alpha(t'')$ and $e = \alpha([t', t''])$.
So, there exist numbers $t_0 = 0 < t_1 < \cdots < t_k = 1$, vertices $v_0=v_x'$, $v_1, \ldots, v_k=v_x''$, and the edges $e_1, \ldots, e_k$ of the tree $\Psi_j$ such that $v_i = \alpha(t_i)$, $i \in \{0, \ldots, k\}$, and $e_i = \alpha([t_{i-1}, t_i])$, $i \in \{1, \ldots, k\}$.

From the choice of the numbers $t'$ and $t''$ it follows that only the points $\alpha(t')$ and $\alpha(t'')$ belong to the intersection of the set $\alpha([t', t''])$ and the subgraph $Q$. Thus $\alpha(t') = v_r$ and $\alpha(t'') = v_s$ for some $r, s \in \{0, \ldots, k\}$, $r < s$. Hence the path $P(v_r, v_s) = \alpha([t', t''])$ connects the vertices $v_r \neq v_s$ of the cycle $Q$ and intersects $Q$ along the set $\{v_r, v_s\}$.

We know already that $\psi(v_r)$, $\psi(v_s) \in \nu_2 \cup \nu_4$. Observe that the points $\psi(v_r)$ and $\psi(v_s)$ can not belong to the different arcs $\nu_2$, $\nu_4$. Really from Proposition~\ref{prop_2.1} it follows that there is $c \in \rr$ such that $\varphi(\Psi_j) \subseteq f^{-1}(c)$, therefore $f \circ \varphi(v_r) = f \circ \varphi(v_s) = c$. On the other hand, as we checked above, the sets $\gamma_2 = \varphi \circ \psi^{-1}(\nu_2)$ and $\gamma_4 = \varphi \circ \psi^{-1}(\nu_4)$ belong to the different level sets of $f$.

Without loss of generality, suppose that $\psi(v_r)$, $\psi(v_s) \in \nu_2$. The set $\nu_2$ is connected, moreover $\nu_2 \subseteq \psi(\Psi_j)$, $\nu_4 \cap \psi(\Psi_j) = \emptyset$ and $\nu_2 \cup \nu_4 \supseteq \psi(Q) \cap \psi (\bigcup_i \Psi_i)$. Therefore the connected set $\psi^{-1}(\nu_2) = Q \cap \Psi_j$ is a subgraph of $P(f)$. Hence there exists the path $\hat{P}(v_r, v_s)$ connecting the vertices $v_r$ and $v_s$ in $Q \cap \Psi_j$.

From the construction we have $P(v_r, v_s) \cup \hat{P}(v_r, v_s) \subseteq \Psi_j$ and $P(v_r, v_s) \neq \hat{P}(v_r, v_s)$. Since $\Psi_j$ is a tree then the vertices $v_r$ and $v_s$ can be connected by the unique path in $\Psi_j$. So, we obtained the contradiction which proves that $\psi(P(f)) \cap \Sigma = \emptyset$. Considering that $\Fr{\Sigma} = \psi(Q) \subseteq \psi(P(f))$ it follows that $\Sigma$ is a connected component of the set $D^2 \setminus \psi(P(f))$.

Suppose now that for some simple cycle $Q' \subseteq P(f)$ the curve $\psi(Q')$ bounds a connected component $\Sigma'$ of the set $D^2 \setminus \psi(P(f))$, but the curve $\varphi(Q')$ is not a boundary of the connected component of the set $D^2 \setminus \varphi(P(f)) = D^2 \setminus P_f$.

Let us prove that in this case $Q' \subseteq \bigcup_i \Psi_i$.

If it does not hold true, then there exists an edge $e_0 \subseteq Q' \cap q(f)$. Evidently,  there exists a connected component $\Theta$ of the set  $D^2 \setminus \varphi(P(f))$ whose boundary contains the set $\varphi(e_0)$. Suppose that $Q = \varphi^{-1}(\Fr{\Theta})$. From Propositions~\ref{prop_2.1} and~\ref{prop_2.2} it follows that $Q$ is a simple cycle. As we proved above the set $\psi(Q)$ is a boundary of some connected component $\Sigma$ of the set $D^2 \setminus \psi(P(f))$. Obviously, $e_0 \subseteq Q \cap Q'$.

Let $x$ be an inner point of an edge $e_0$, $z = \psi(x)$. By the conditions of proposition we have $z \in \psi(q(f)) = \partial D^2$. It is easy to see that for sufficiently small neighborhood $W$ of the point $z$ in $D^2$ which is homeomorphic to half-disk the set $W \setminus \psi(P(f)) = W \setminus \psi(e_0)$ is connected. Therefore $W \setminus \psi(P(f)) \subseteq \Sigma \cap \Sigma' \neq \emptyset$. From $\Sigma \cap \Fr{\Sigma'} \subseteq \Sigma \cap \psi(P(f)) = \emptyset$ it follows that $\Sigma \subseteq \Sigma'$. By a parallel argument $\Sigma' \subseteq D^2 \setminus \psi(P(f))$, thus $\Sigma' \cap \Fr{\Sigma} \subseteq \Sigma' \cap \psi(P(f)) = \emptyset$ and $\Sigma' \subseteq \Sigma$.
Hence $\Sigma' = \Sigma$, $Q' = Q$ and the curve $\varphi(Q')$ bounds a connected component of the set $D^2 \setminus P_f$, but it contradicts to the choice of the cycle $Q'$. Therefore $Q' \subseteq \bigcup_i \Psi_i$.

The set $Q'$ is connected thus there is $j$ such that $Q' \subseteq \Psi_j$. But $\Psi_j$ is tree and no one cycle is contained in it. This contradiction is a final step of proof.
\end{proof}

\begin{nas}
In the conditions of Lemma~\ref{lemma_2_2} there exists a homeomorphism $\Phi : D^{2} \rightarrow D^{2}$ such that
$\Phi \circ \psi_{1} = \psi_{2}$.
\end{nas}
\begin{proof}
Let $\Sigma_1, \ldots, \Sigma_k$ be the connected components of the set $D^2 \setminus \psi_1(P(f))$ and $Q_1, \ldots, Q_k$ be the cycles of the graph $P(f)$ such that $\psi_1(Q_i) = \Fr{\Sigma_i}$, $i \in \{1, \ldots, k\}$, see Proposition~\ref{prop_2.2}. It is clear that $P(f) = \bigcup_{i=1}^k Q_i$.

By Lemma~\ref{lemma_2_2}, every set $\psi_2(Q_i)$ is a boundary of some connected component $\Sigma_i'$ of $D^2 \setminus \psi_2(P(f))$. By using Lemma~\ref{lemma_2_2} once again it is easy to see that $D^2 \setminus \psi_2(P(f)) = \bigcup_{i=1}^k \Sigma_i'$.
By Schoenflies's theorem, for every $i \in \{1, \ldots, k\}$,  the homeomorphism $\psi_2 \circ\ \psi_1^{-1}\restrict{\vphantom{\psi_1}}{\psi_1(Q_i)} : \psi_1(Q_i) \rightarrow \psi_2(Q_i)$ can be extended to a homeomorphism of disks $\Phi_i : \Cl{\Sigma_i} \rightarrow \Cl{\Sigma_i'}$, see ~\cite{Newman}. It is easy to see that the map $\Phi : D^2 \rightarrow D^2$,
\[
\Phi(z) = \Phi_i(z)\,, \quad \mbox{for } z \in \Cl{\Sigma_i} \,,
\]
is well defined and maps $D^2$ onto itself bijectively. A finite family of the closed sets $\{\Cl{\Sigma_i}\}$ generates the fundamental cover of $D^2$ and on each of them $\Phi$ is continuous. Hence $\Phi$ is continuous on $D^2$, see~\cite{RF}. It is known that a continuous bijective map of compactum to a Hausdorff space is a homeomorphism.
\end{proof}

\section{The conditions of topological equivalence}
Let $f : D^2 \rightarrow \rr$ be a pseudoharmonic function and $a_1 < \cdots < a_N$ be all its critical and semiregular values. Let us consider a homeomorphism $h_f : [a_1, a_N] \rightarrow [1, N]$ such that $h_f(a_j) = j$ for all $j \in \{1, \ldots, N\}$. It is easy to see that a continuous function $\hat{f} = h_f \circ f$ is pseudoharmonic, set of its critical and semiregular values is $\{1, \ldots, N\}$, and $P(\hat{f}) = P(f)$. Function $\hat{f}$ is called \emph{a standardization} of $f$.

\begin{theorem}\label{main_th} Two pseudoharmonic functions $f$ and $g$  are topologically equivalent iff there exists an isomorphism of combinatorial diagrams $\varphi:P(f)\rightarrow P(g)$ which preserves a strict partial order defined on them and the orientation.
\end{theorem}
\begin{proof}
\textit{Necessity.}
Suppose that two pseudoharmonic functions $f: D^2 \rightarrow \mathbb{R}$ and $g: D^2 \rightarrow \mathbb{R}$ are topologically equivalent. Then there exist homeomorphisms $H : D^2 \rightarrow D^2$ and $h : \mathbb{R}\rightarrow \mathbb{R}$ such that $f = h^{-1} \circ\ g \circ H$. Also to $f$ and $g$ correspond their combinatorial diagrams $P(f)$ and $P(g)$ with the strict partial order and the orientation which conform to $f$ and $g$. Let $\psi_{1}$ and $\psi_{2}$ be embeddings of $P(f)$ and $P(g)$ into $D^{2}$ which are consistent with $f$ and $g$, respectively (recall that the partial orientations on $Cr$-subgraphs of $P(f)$ and $P(g)$ are the same as the orientation of $\partial D^2$).
Evidently, the homeomorphism $H$ maps the regular points of $f$ onto the regular points of $g$ and the critical points of $f$ onto the critical points of $g$, respectively. From $h \circ f = g \circ H$ and bijectivity of $h$ it follows that the homeomorphism $H$ maps the regular, critical and semiregular levels of $f$ onto regular, critical and semiregular levels of $g$, respectively. Thus $H \circ \psi_{1}(P(f)) = \psi_{2}(P(g))$ and the bijective map $\varphi = \psi_{2}^{-1} \circ H \circ \psi_{1} : P(f) \rightarrow P(g)$ is defined. So, we have the following commutative diagram:
\[
\begin{CD}
P(f)	@>{\psi_{1}}>>	D^{2}	@>{f}>>	\rr \\
@V{\varphi}VV		@V{H}VV	@VV{h}V \\
P(g)	@>>{\psi_{2}}>	D^{2}	@>>{g}>	\rr
\end{CD}
\]
It is easy to see that $\varphi$ defines an isomorphism of graphs. Let us prove that the maps $\varphi$ and $\varphi^{-1}$ are monotone.
We should remind that only the preserving orientation homeomorphisms $\rr$  are considered thus the map $h : \rr \rightarrow \rr$ preserves an order of points of $\rr$. Let $v_{1}$ and $v_{2}$ be two vertices of the diagram $P(f)$. By the definition of the diagram $P(f)$ an inequality $v_{1} < v_{2}$ is equivalent to $f \circ \psi_{1}(v_{1}) < f \circ \psi_{1}(v_{2})$, so, that is also equivalent to $h \circ f \circ \psi_{1}(v_{1}) < h \circ f \circ \psi_{1}(v_{2})$. This inequality is equivalent to $g \circ \psi_{2} \circ \varphi(v_{1}) < g \circ \psi_{2} \circ \varphi(v_{2})$ since $h \circ f \circ \psi_{1} = g \circ H  \circ \psi_{1} = g \circ \psi_{2} \circ \varphi$. By the definition of the relation of order on $P(g)$, the last inequality is equivalent to $\varphi(v_{1}) < \varphi(v_{2})$. So, the inequalities $v_{1} < v_{2}$ and $\varphi(v_{1}) < \varphi(v_{2})$ are equivalent. Finally, we remind that  $\varphi$ is bijective by the construction thus $\varphi$ and $\varphi^{-1}$ are monotone. From what we said it follows that $\varphi:P(f)\rightarrow P(g)$ is an isomorphism of diagrams.

\medskip

\textit{Sufficiency.} Suppose that $f, g: D^2 \rightarrow \mathbb{R}$ are pseudoharmonic functions and $P(f)$, $P(g)$ are their diagrams such that there exists an isomorphism $\phi:P(f)\rightarrow P(g)$ preserving a strict partial order on the set of vertices and the partial orientations on their $Cr$-- subgraphs.

At first, we want to replace the functions $f$ and $g$ on the normalized pseudoharmonic functions $\hat{f}$ and $\hat{g}$ with the same combinatorial diagrams as $f$ and $g$.

Let $a_1 < \cdots < a_N$ be the critical and the semiregular values of $f$ and $b_1 < \cdots < b_M$ be critical and semiregular values of $g$. We fix homeomorphisms $h_f : [a_1, a_N] \rightarrow [1, N]$ and $h_g : [b_1, b_M] \rightarrow [1, M]$ such that $h_f(a_i) = i$ and $h_g(b_j) = j$ for all $i \in \{1, \ldots, N\}$ and $j \in \{1, \ldots, M\}$. Obviously, the maps $h_f$ and $h_g$ are orientation preserving. Required normalized pseudoharmonic functions have the following forms $\hat{f} = h_f \circ f$ and $\hat{g} = h_g \circ g$.

Let us fix an embedding $\psi_f : P(f) \rightarrow D^2$ which is consistent with $f$ and an embedding $\psi_g : P(g) \rightarrow D^2$ which is consistent with $g$ (we should remark that the embeddings $\psi_f$ and $\psi_g$ are also consistent with $\hat{f}$ and $\hat{g}$, respectively).

Let us prove that for any vertex $v$ of $P(f)$ the following condition holds true
\begin{equation}\label{eq_14}
\hat{f} \circ \psi_f(v) = \hat{g} \circ \psi_g \circ \phi(v) \,.
\end{equation}
We remark that
\[
f \circ \psi_f(V(P(f))) = \{a_1, \ldots, a_N\} \,, \quad g \circ \psi_g(V(P(g))) = \{b_1, \ldots, b_M\} \,.
\]
Hence
\[
\hat{f} \circ \psi_f(V(P(f))) = \{1, \ldots, N\} \,, \quad \hat{g} \circ \psi_g(V(P(g))) = \{1, \ldots, M\} \,.
\]
Fix the sequence of vertices $u_1, \ldots, u_s \in V(P(f))$ such that $\hat{f} \circ \psi_f(u_1) = 1, \ldots, \hat{f} \circ \psi_f(u_s) = s = \hat{f} \circ \psi_f(v)$. Then $u_1 < \cdots < u_s$ in $P(f)$ hence $\phi(u_1) < \cdots < \phi(u_s)$ in $P(g)$ and $\hat{g} \circ \psi_g(u_1) \circ \phi < \cdots < \hat{g} \circ \psi_g \circ \phi(u_s)$. Thus $j \leq \hat{g} \circ \psi_g \circ \phi(u_j)$, $j \in \{1, \ldots, s\}$. Hence $\hat{f} \circ \psi_f(v) = \hat{f} \circ \psi_f(u_s) = s \leq \hat{g} \circ \psi_g \circ \phi(u_s) = \hat{g} \circ \psi_g \circ \phi(v)$.
By replacing $f$ at $g$, we have $\hat{f} \circ \psi_f(v) \geq \hat{g} \circ \psi_g \circ \phi(v)$. So,~\eqref{eq_14} holds true.
From~\eqref{eq_14} it follows  that $M = N$ and $\hat{f}(D^2) = \hat{g}(D^2) = [1, N]$.

Let us construct a homeomorphism $\varphi : P(f) \rightarrow P(g)$ such that it realizes an isomorphism $\phi$ and satisfies the following relation on the space $P(f)$
\begin{equation}\label{eq_15}
\hat{f} \circ \psi_f = \hat{g} \circ \psi_g \circ \varphi \,.
\end{equation}
By the definition we have $\varphi(v) = \phi(v)$, $v \in V(P(f))$, on the set of vertices.

Suppose that an edge $e = e(v', v'') \in E(P(f))$ belongs to a subgraph $\bigcup_i \Psi_i(f) = \Cl{P(f) \setminus q(f)}$. The set $e$ is connected and $\hat{f} \circ \psi_f$ is locally constant on $\bigcup_i \Psi_i(f)$, therefore $\psi_f(e) \subseteq \hat{f}^{-1}(c)$ for some $c \in \rr$. In particular, $\hat{f} \circ \psi_f(v') = \hat{f} \circ \psi_f(v'') = c$ and vertices $v'$ and $v''$ are non comparable in $P(f)$. So, the vertices $\phi(v')$ and $\phi(v'')$ are non comparable in $P(g)$ and $\phi(e) \subset \bigcup_j \Psi_j(g) = \Cl{P(g) \setminus q(g)}$. Hence $\psi_g \circ \phi(e) \subset \hat{g}^{-1}(c')$ for some $c' \in \rr$, in particular, $\hat{g} \circ \psi_f \circ \phi(v') = \hat{g} \circ \psi_f \circ \phi(v'') = c'$. But from~\eqref{eq_14} it follows that $\hat{g} \circ \psi_f \circ \phi(v') = \hat{f} \circ \psi_f(v') = c$, therefore $e \subset (\hat{f} \circ \psi_f)^{-1}(c)$ and $\phi(e) \subset (\hat{g} \circ \psi_g)^{-1}(c)$.
Fix a homeomorphism $\varphi_e : e \rightarrow \phi(e)$ such that $\varphi(v') = \phi(v')$ and $\varphi(v'') = \phi(v'')$. Obviously, $\hat{f} \circ \psi_f(x) = \hat{g} \circ \psi_g \circ \varphi_e(x) = c$, $x \in e$.

Suppose that an edge $e = e(v', v'') \in E(P(f))$ belongs to $\mathcal{C}r$-subgraph $q(f)$. Then every point of a set $\psi_f(e) \setminus \{\psi_f(v'), \psi_f(v'')\}$ is a regular boundary point of $\hat{f}$, hence $\hat{f}$ is strictly monotone on the arc $\psi_f(e)$ and maps it homeomorphically on $[c', c'']$, where
\[
c' = \min(\hat{f} \circ \psi_f(v'), \hat{f} \circ \psi_f(v'')) \,, \quad
c'' = \max(\hat{f} \circ \psi_f(v'), \hat{f} \circ \psi_f(v'')) \,.
\]
Since $\phi$ is an isomorphism of the combinatorial diagrams then from $C1$ it follows that $\phi(e) \in q(g)$. Thus $\hat{g}$ maps the set $\psi_g(\phi(e))$ in $\rr$ homeomorphically. From~\eqref{eq_14} it follows that $\hat{g} \circ \psi_g(\phi(e)) = [c', c'']$.

Suppose that
\[
\varphi_e = \Bigl(\hat{g} \circ \psi_g\restrict{\vphantom{g}}{\phi(e)}\Bigr)^{-1} \circ \hat{f} \circ \psi_f: e \rightarrow \phi(e) \,.
\]
It is easy to see that this map is a homeomorphism and satisfies the following relation $\hat{f} \circ \psi_f(x) = \hat{g} \circ \psi_g \circ \varphi_e(x)$, $x \in e$.

Let us define a map $\varphi : P(f) \rightarrow P(g)$ as
\[
\varphi(x) = \varphi_e(x) \,, \quad \mbox{for } x \in e \,.
\]
By the construction $\varphi_e(v) = \phi(v)$ for $v \in e \cap V(P(f))$ therefore $\varphi_{e'}(x) = \varphi_{e''}(x)$ for every pair of edges $e', e'' \in E(P(f))$ and $x \in e' \cap e'' \subseteq V(P(f))$. So, the map $\varphi$ is defined correctly. It is easy to see that $\varphi$ satisfies~\eqref{eq_15}.
The collection of edges $\{ e \in E(P(f)) \}$ generate a finite closed cover of a space $P(f)$ thus it is fundamental. Hence $\varphi$ is continuous since each of maps $\varphi_e$ is continuous by definition, where $e \in E(P(f))$, see~\cite{RF}.
It is easy to see that $\varphi$ is a bijective map and the spaces $P(f)$ and $P(g)$ are compact. Therefore $\varphi$ maps $P(f)$ on $P(g)$ homeomorphically. Moreover, since $\phi$ preserves orientation of $q(f)$, then an orientation on $q(g)=\varphi(q(f))$ induced by $\varphi$ coincides with the orientation of $q(g)$ in $P(g)$.

We set
\[
H_0 = \psi_g \circ \varphi \circ \psi_f^{-1}\restrict{\vphantom{f}}{P_f} : P_f \rightarrow P_g \,.
\]
By the construction $H_0$ maps the set $P_f = \psi_f(P(f))$ on $P_g = \psi_g(P(g))$ homeomorphically. Moreover, from ~\eqref{eq_15} it follows that
\begin{equation}\label{eq_16}
\hat{g} \circ H_0 = \hat{f} \,.
\end{equation}

As orientations induced on $\partial D^2=\psi_{f}(q(f))=\psi_{g}(q(g))$ by $\psi_f$ and $\psi_g$ from $q(f)$ and $q(g)$ respectively coincide with the positive orientation of $\partial D^2$ by definition, then $H_0$ preserves the orientation of $\partial D^2$.

Our aim is to extend $H_0$ to a homeomorphism $H : D^2 \rightarrow D^2$ such that $\hat{g} \circ H = \hat{f}$.

Let $\Theta$ be one of connected components of $D^2 \setminus P(f)$. From Propositions~\ref{prop_2.1} and~\ref{prop_2.2} it follows that there exists a simple cycle $Q \subseteq P(f)$ such that $\psi_f(Q) = \Fr{\Theta}$.
Evidently, $\phi(q(f)) = \varphi(q(f)) = q(g)$. Thus $\psi_g \circ \varphi(q(f)) = \psi_g(q(g))= \partial D^2$ and from Lemma~\ref{lemma_2_2} it follows that a set $\psi_g \circ \varphi(Q)$ is a boundary of some connected component $\Sigma$ of $D^2 \setminus P_g$.

Denote by $\Rfr(f)$ a set of all regular boundary points of $\hat{f}$. It is easy to see that, on one hand, $\Rfr(f) = \psi_f(q(f) \setminus V(P(f)))$, on the other hand the set $W_f$ of all regular boundary points of $\hat{f}\restrict{\vphantom{f}}{\Cl{\Theta}}$ coincides with $\Rfr(f) \cap \Cl{\Theta} = \Rfr(f) \cap \Fr{\Theta} = \Rfr(f) \cap \psi_f(Q)$. Therefore $W_f$ is an image of a set $Q \cap (q(f) \setminus V(P(f)))$.

Let $\Rfr(g)$ be a set of regular boundary points of $\hat{g}$. By analogy, we can conclude that the set $W_g$ of regular boundary points of $\hat{g}\restrict{\vphantom{g}}{\Cl{\Sigma}}$ is an image of a set $\varphi(Q) \cap (q(g) \setminus V(P(g)))$. But a map $\varphi$ is bijective and, also, it is known that $q(g) = \varphi(q(f))$ and $V(P(g)) = \varphi(V(P(f)))$. Therefore $\varphi(Q) \cap (q(g) \setminus V(P(g))) = \varphi(Q \cap (q(f) \setminus V(P(f))))$ and $W_g = \psi_g \circ \varphi \circ \psi_f^{-1}(W_f) = H_0(W_f)$.

From Lemma~\ref{lemma_2_C4} it follows that the function $\hat{f}$ is regular on the set $\Cl{\Theta}$. Let $z_1, \ldots, z_4$ and $\gamma_1, \ldots, \gamma_4$ be the points and the arcs, respectively, from Definition~\ref{ozn_weak_regular}. Proposition~\ref{prop_min_num_of_arcs} guarantees that $W_f = \mathring{\gamma}_1 \cup \mathring{\gamma}_3$ holds true. We set $K_f = \gamma_2 \cup \gamma_4 = \Fr{\Theta} \setminus W_f$.

Similarly, the function $\hat{g}$ is regular on the set $\Cl{\Sigma}$. Let $w_1, \ldots, w_4$ and $\nu_1, \ldots, \nu_4$ be the points and the arcs, respectively, from Definition~\ref{ozn_weak_regular}. Then $W_g = \mathring{\nu}_1 \cup \mathring{\nu}_3$. We set $K_g = \nu_2 \cup \nu_4 = \Fr{\Sigma} \setminus W_g$.

We already verified that $W_g = H_0(W_f)$. The map $H_0$ is bijective, thus $K_g = H_0(K_f)$. Hence for the functions $\hat{f}\restrict{\vphantom{f}}{\Cl{\Theta}}$ and $\hat{g}\restrict{\vphantom{g}}{\Cl{\Sigma}}$ Theorem~\ref{theorem_1.3} is satisfied with the same set $D' \in \{I^2, \Cl{D}^{2}_{+}, D^2\}$ and its subset
\begin{gather*}
K' = \bigl\{ (x, y) \in D' \,|\, y \in \{y_1, y_2\} \bigr\} \,;\\
y_1 = \min \{ y \,|\, (x, y) \in D' \} \,,\\
y_2 = \max \{ y \,|\, (x, y) \in D' \} \,.
\end{gather*}
Fix a homeomorphism $\chi_f : \Fr{\Theta} \rightarrow \Fr{D'}$ such that $\chi_f(K_f) = K'$. We set
\[
\chi_g = \chi_f \circ \psi_f \circ \varphi^{-1} \circ \psi_g^{-1} = \chi_f \circ H_0^{-1} : \Fr{\Sigma} \rightarrow \Fr{D'} \,.
\]
The map $\chi_g$ is a composition of homeomorphisms therefore $\chi_g$ is a homeomorphism. Moreover, $\chi_g(K_g) = \chi_f \circ H_0^{-1}(K_g) = \chi_f(K_f) = K'$.

From Theorem~\ref{theorem_1.3} it follows that there exist numbers $a_f, b_f, a_g, b_g \in \rr$ and homeomorphisms $F_Q : \Cl{\Theta} \rightarrow D'$ and $G_Q : \Cl{\Sigma} \rightarrow D'$ such that $F_Q\restrict{\vphantom{F}}{K_f} = \chi_f$, $G_Q\restrict{\vphantom{G}}{K_g} = \chi_g$, and $\hat{f} \circ F_Q^{-1}(x, y) = a_f y + b_f$, $\hat{g} \circ G_Q^{-1}(x, y) = a_g y + b_g$, $(x, y) \in D'$.
We set
\begin{gather*}
K_1 = \{ (x, y) \in D' \,|\, y = y_1 \} \,,\\
K_2 = \{ (x, y) \in D' \,|\, y = y_2 \} \,.
\end{gather*}
Due to the choice of $D'$ the sets $K_1$ and $K_2$ are connected. Hence $F_Q^{-1}(K_i)$, $i = 1, 2$, are also connected.
From $F_Q^{-1}(K_1) = \chi_f^{-1}(K_1) \subseteq K_f \subseteq \psi_f(V(P(f)) \cup \bigcup_i \Psi_i(f))$ it follows that there exists $c_1 \in \rr$ such that $F_Q^{-1}(K_1) \subseteq \hat{f}^{-1}(c_1)$ (we remind that $\hat{f}$ is locally constant on the set $\psi_f(V(P(f)) \cup \bigcup_i \Psi_i(f))$). On the other hand, $G_Q^{-1}(K_1) = \chi_g^{-1}(K_1) = (\chi_f \circ H_0^{-1})^{-1}(K_1)$, therefore $\hat{g} \circ G_Q^{-1}(K_1) = \hat{g} \circ H_0 \circ \chi_f^{-1}(K_1) = \hat{f} \circ \chi_f^{-1}(K_1) = \hat{f} \circ F_Q^{-1}(K_1) = c_1$ since~\eqref{eq_16} and $G_Q^{-1}(K_1) \subseteq \hat{g}^{-1}(c_1)$.
Similarly, there exists $c_2 \in \rr$ such that $F_Q^{-1}(K_2) \subseteq \hat{f}^{-1}(c_2)$ and $G_Q^{-1}(K_2) \subseteq \hat{g}^{-1}(c_2)$.
Hence for any $(x_1, y_1) \in K_1$ and $(x_2, y_2) \in K_2$ the following conditions hold true
\begin{equation}\label{eq_17}
\left\{
\begin{aligned}
\hat{f} \circ F_Q^{-1}(x_1, y_1) & = a_f y_1 + b_f = c_1 \,,\\
\hat{g} \circ G_Q^{-1}(x_1, y_1) & = a_g y_1 + b_g = c_1 \,,\\
\hat{f} \circ F_Q^{-1}(x_2, y_2) & = a_f y_2 + b_f = c_2 \,,\\
\hat{g} \circ G_Q^{-1}(x_2, y_2) & = a_g y_2 + b_g = c_2 \,.
\end{aligned}
\right.
\end{equation}
It is easy to see that a determinant of this system of linear equations with variables $a_f$, $b_f$, $a_g$ and $b_g$ equals to $(y_2 - y_1)^2$. By the construction $y_1 \neq y_2$ thus $(y_2 - y_1)^2 \neq 0$ and the system~\eqref{eq_17} has the unique solution which can be easily calculated:
\[
a_f = a_g = \frac{c_2-c_1}{y_2-y_1} \,, \quad b_f = b_g = \frac{c_1 y_2-c_2 y_1}{y_2-y_1} \,.
\]
So, on the set $D'$ the following equality holds true
\begin{equation}\label{eq_18}
\hat{f} \circ F_Q^{-1} = \hat{g} \circ G_Q^{-1} \,.
\end{equation}
It is clear that $G_Q^{-1} \circ F_Q(\Fr{\Theta}) = \Fr{\Sigma}$. We remind that
\[
G_Q^{-1} \circ F_Q\restrict{\vphantom{F}}{K_f} = \chi_g^{-1} \circ \chi_f\restrict{\vphantom{a}}{K_f} = H_0\restrict{\vphantom{H}}{K_f} \,.
\]
Since $G_Q^{-1} \circ F_Q(K_f) = H_0(K_f) = K_g$, then a homeomorphism $G_Q^{-1} \circ F_Q$ satisfies to relations $G_Q^{-1} \circ F_Q(W_f) = G_Q^{-1} \circ F_Q(\Fr{\Theta} \setminus K_f) = \Fr{\Sigma} \setminus K_g = W_g$.

As we know, the set $W_f$ has two connected components $\mathring{\gamma}_1$ and $\mathring{\gamma}_3$. Under the action of the homeomorphism $H_0$ they have to map on the connected components $\mathring{\nu}_1$ and $\mathring{\nu}_3$ of the set $W_g$. We cyclically change a numeration of the points $w_1, \ldots, w_4$ and the arcs $\nu_1, \ldots, \nu_4$ so that $H_0(\mathring{\gamma}_{2k-1}) = \mathring{\nu}_{2k-1}$, $k = 1, 2$.
The homeomorphism $G_Q^{-1} \circ F_Q$ also has to map the sets $\mathring{\gamma}_1$ and $\mathring{\gamma}_3$ on the connected components of the set $W_g$.

Let us prove that under the condition $\mathring{\gamma}_2 \cup \mathring{\gamma}_4 \neq \emptyset$ the relations hold true $G_Q^{-1} \circ F_Q(\mathring{\gamma}_{2k-1}) = \mathring{\nu}_{2k-1}$, $k = 1, 2$.
We remark that either $G_Q^{-1} \circ F_Q(\mathring{\gamma}_1) = \mathring{\nu}_1$ or $G_Q^{-1} \circ F_Q(\mathring{\gamma}_1) = \mathring{\nu}_3$ holds true.

It is clear that an arc $\gamma_{2k-1}$ is a closure of a arc $\mathring{\gamma}_{2k-1}$ in $D^2$ (by definition $\mathring{\gamma}_{2k-1} \neq \emptyset$), $k = 1, 2$; similarly, $\nu_{2k-1} = \Cl{\mathring\nu}_{2k-1}$. Thus if $G_Q^{-1} \circ F_Q(\gamma_{2k-1}) \neq \nu_{2j-1}$ for some $k, j \in \{1, 2\}$, then  $G_Q^{-1} \circ F_Q(\mathring\gamma_{2k-1}) \neq \mathring\nu_{2j-1}$.

Without loss of generality, suppose that $\mathring\gamma_2 \neq \emptyset$. Then $z_2 \in \gamma_1 \setminus \gamma_3$ and $H_0(z_2) \in \nu_1 \setminus \nu_3$. But $z_2 \in \gamma_2 \subset K_f$ and $H_0(z_2) = G_Q^{-1} \circ F_Q(z_2)$. Thus $G_Q^{-1} \circ F_Q(\gamma_1) \neq \nu_3$, hence $G_Q^{-1} \circ F_Q(\mathring\gamma_1) = \mathring\nu_1$ and $G_Q^{-1} \circ F_Q(\mathring\gamma_3) = \mathring\nu_3$.

Suppose now that $\mathring{\gamma}_2 \cup \mathring{\gamma}_4 = \emptyset$. Then both sets $\gamma_2$ and $\gamma_4$ are one-point and $D' = D^2$. Consider an involution $Inv : D^2 \rightarrow D^2$, $Inv(x, y) = (-x, y)$, $(x, y) \in D^2$. Obviously, it changes the connected components of $\Fr{D'} \setminus K' = F_Q(W_f)$. Moreover, $\restrict{Inv}{K'} = Id$ since $K' = \{(0, -1), (0, 1)\}$.
If $G_Q^{-1} \circ F_Q(\mathring\gamma_1) = \nu_3$, then the map $G_Q$ can be replaced by $Inv \circ G_Q$. It is easy to see that the following conditions hold true
\begin{itemize}
	\item $Inv \circ G_Q\restrict{\vphantom{G}}{K_g} = \chi_g\restrict{\vphantom{a}}{K_g}$;
	\item $\hat{g} \circ (Inv \circ G_Q)^{-1}(x, y) = \hat{g} \circ G_Q^{-1}(-x, y) = \hat{g} \circ G_Q^{-1}(x, y) = a_g y + b_g$;
	\item $(Inv \circ G_Q)^{-1} \circ F_Q(\mathring\gamma_1) = W_g \setminus \mathring\nu_3 = \mathring\nu_1$.
\end{itemize}

So, we proved that the homeomorphisms $F_Q : \Cl{\Theta} \rightarrow D'$ and $G_Q : \Cl{\Sigma} \rightarrow D'$ satisfy conditions
\begin{itemize}
	\item $\hat{f} \circ F_Q^{-1} = \hat{g} \circ G_Q^{-1}$;
	\item $G_Q^{-1} \circ F_Q\restrict{\vphantom{F}}{K_f} = H_0\restrict{\vphantom{H}}{K_f}$;
	\item $G_Q^{-1} \circ F_Q(\mathring{\gamma}_{2k-1}) = \mathring{\nu}_{2k-1}$, $k = 1, 2$.
\end{itemize}
We set $H_Q = G_Q^{-1} \circ F_Q : \Cl{\Theta} \rightarrow \Cl{\Sigma}$. Let us verify that
\[
\restrict{H_Q}{\Fr{\Theta}} = H_0\restrict{\vphantom{H_Q}}{\Fr{\Theta}}\,.
\]
It is sufficient to prove that $H_Q(z) = H_0(z)$ for all $z \in W_f = \Fr{\Theta} \setminus K_f \subset \gamma_1 \cup \gamma_3$.
As we know, the set $\mathring\gamma_1$ consists of the regular boundary points of the function $\hat{f}$, therefore $\hat{f}$ is strictly monotone on the arc $\gamma_1$ and maps it on $\hat{f}(\gamma_1) \subset \rr$ homeomorphically ( since $\gamma_1$ is the compactum and the space $\rr$ is Hausdorff). Similarly, a map $\hat{g}\restrict{\vphantom{g}}{\nu_1} : \nu_1 \rightarrow \hat{g}(\nu_1) \subset \rr$ is a homeomorphism onto its image.

As a consequence of $\gamma_1 \subset P_f$, $\nu_1 \subset P_g$ and from~\eqref{eq_16} it follows that $\hat{f}(\gamma_1) = \hat{g} \circ H_0(\gamma_1) = \hat{g}(\nu_1)$. Thus the following map is well defined
\[
\hat{g}^{-1} \circ \hat{f}\restrict{\vphantom{f}}{\gamma_1} : \gamma_1 \rightarrow \nu_1 \,.
\]
By using~\eqref{eq_16} again we have $H_0\restrict{\vphantom{H}}{\gamma_1} = \hat{g}^{-1} \circ \hat{f} \restrict{\vphantom{f}}{\gamma_1}$.
On the other hand, from~\eqref{eq_18} it follows that $\hat{f} = \hat{g} \circ G_Q^{-1} \circ F_Q = \hat{g} \circ H_Q$, therefore $\hat{g}^{-1} \circ \hat{f}\restrict{\vphantom{f}}{\gamma_1} = H_Q\restrict{\vphantom{H}}{\gamma_1}$. Hence $H_0\restrict{\vphantom{H}}{\gamma_1} = H_Q\restrict{\vphantom{H}}{\gamma_1}$.

By analogy we prove that $H_0\restrict{\vphantom{H}}{\gamma_3} = H_Q\restrict{\vphantom{H}}{\gamma_3}$.
So, we constructed the homeomorphism $H_Q : \Cl{\Theta} \rightarrow \Cl{\Sigma}$ such that
\[
\hat{f} = \hat{g} \circ H_Q
\]
and $H_Q\restrict{\vphantom{H}}{\Fr{\Theta}} = H_0\restrict{\vphantom{H}}{\Fr{\Theta}}$.

Let us construct a homeomorphism $H : D^2 \rightarrow D^2$ such that $\hat{f} = \hat{g} \circ H$.
For every connected component $\Theta$ of $D^2 \setminus P_f$ its boundary $\Fr{\Theta}$ is an image of a simple cycle $Q(\Theta)$, see Propositions ~\ref{prop_2.1} and~\ref{prop_2.2}. For every $\Theta$ we fix the homeomorphism $H_{Q(\Theta)}$ such that $\hat{f} = \hat{g} \circ H_{Q(\Theta)}$ and $H_{Q(\Theta)}\restrict{\vphantom{H}}{\Fr{\Theta}} = H_0\restrict{\vphantom{H}}{\Fr{\Theta}}$.
We define $H$ by the following relations
\[
H(z) = H_{Q(\Theta)} \,, \quad \mbox{for } z \in \Cl{\Theta} \,.
\]
If $z \in \Cl{\Theta'} \cap \Cl{\Theta''}$, then $z \in P_f$ and $H_{Q(\Theta')}(z) = H_0(z) = H_{Q(\Theta'')}(z)$. So, the map $H$ is correctly defined.
By Lemma~\ref{lemma_2_2} there exists a bijective correspondence between the connected components of the sets $D^2 \setminus P_f$ and $D^2 \setminus P_g$. Thus $H(\Theta') \cap H(\Theta'') = \emptyset$ for $\Theta' \neq \Theta''$ and $\bigcup_{\Theta} H(\Cl{\Theta}) = D^2$. Hence the map $H$ is bijective.

Evidently, by the construction we have
\[
\hat{f} = \hat{g} \circ H \,.
\]
The closures of the connected components of $D^2 \setminus P_f$ generate a finite closed cover of disk $D^2$. It is known~\cite{RF} that this cover is fundamental. Therefore the map $H$ is continuous on $D^2$ since by construction it is continuous on each element of its cover.

It is known that a continuous bijective map of compactum in Hausdorff's space is a homeomorphism. So, $H : D^2 \rightarrow D^2$ is a homeomorphism.

Now recall that map $H_0\mid_{\partial D^2}=H\mid_{\partial D^2}$ preserves the orientation. Consequently, $H$ preserves the orientation on $D^2$.

We remind that $\hat{f} = h_f \circ f$ and $\hat{g} = h_g \circ g$ for some homeomorphisms $h_f : f(D^2) = [a_1, a_N] \rightarrow [1, N]$ and $h_g : g(D^2) = [b_1, b_N] \rightarrow [1, N]$ which preserve orientation. It is obvious that the map $h_0 = h_g^{-1} \circ h_f : f(D^2) \rightarrow g(D^2)$ is a homeomorphism of the segment $f(D^2)$ on the segment $g(D^2)$ which preserves the orientation. Let us fix a homeomorphism $h : \rr \rightarrow \rr$, which preserves the orientation and satisfies $\restrict{h}{f(D^2)} = h_0$.
It is easy to see that
\[
h \circ f = h_g^{-1} \circ h_f \circ f = h_g^{-1} \circ \hat{f} = h_g^{-1} \circ \hat{g} \circ H = g \circ H \,,
\]
so, the functions $f$ and $g$ are topologically equivalent.
\end{proof}

\begin{figure}[htbp]
\begin{center}
\includegraphics{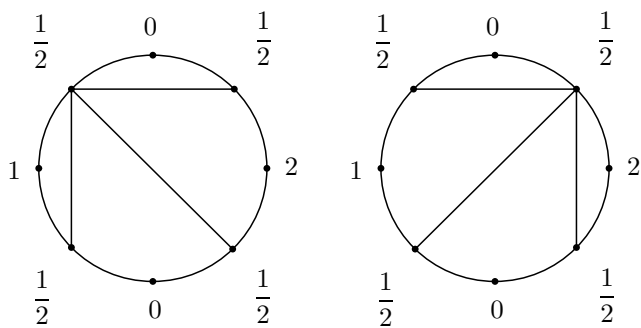}
\caption{The diagrams of topologically non equivalent pseudoharmonic functions.}
\end{center}
\label{fig3}
\end{figure}

On Fig.~\ref{fig3} the diagrams of two pseudoharmonic functions which have two local minima, two local maxima on $\partial D^{2}$ and one boundary critical point are represented. But, these two functions are not topologically equivalent.

\end{document}